\documentclass[preprint,12pt,authoryear]{elsarticle}

\usepackage{placeins}
\usepackage[normalem]{ulem}
\usepackage{amsmath}
\usepackage{setspace}
\usepackage[hidelinks]{hyperref}
\usepackage[utf8]{inputenc}
\usepackage[T1]{fontenc}
\usepackage{hyperref}
\usepackage{url}
\usepackage{booktabs}
\usepackage{amsmath}
\usepackage{array}
\usepackage{amsthm}
\usepackage{amsfonts}
\usepackage{nicefrac}
\usepackage{microtype}
\usepackage{graphicx}
\usepackage{doi}
\usepackage{xcolor}
\usepackage{pifont}
\usepackage[ruled,vlined]{algorithm2e}
\usepackage{tikz}
\usetikzlibrary{positioning,arrows.meta,fit,backgrounds,calc}
\usepackage[most]{tcolorbox}
\usepackage{amssymb}
\usepackage{listings}
\usepackage{float}
\newtheorem{theorem}{Theorem}

\tcbuselibrary{breakable,skins}

\tcbset{
  pmcbox/.style={
    enhanced,
    rounded corners,
    arc=2mm,
    boxrule=0.9pt,
    colback=blue!2,
    colframe=blue!45!black,
    colbacktitle=blue!75!black,
    coltitle=white,
    left=5pt,right=5pt,top=5pt,bottom=5pt,
    fonttitle=\bfseries,
    drop shadow=black!15,
  }
}

\newcommand{\cmark}{\textcolor{green!60!black}{\ding{51}}}
\newcommand{\xmark}{\textcolor{red!70!black}{\ding{55}}}

\journal{European Journal of Operational Research}

\begin{document}
\onehalfspacing
\begin{frontmatter}

\title{Solver-Verified Formulation Generation and Selection for Multi-Warehouse Inventory Allocation Using Large Language Models}

\author[inst1]{Jintao Xu}
\ead{xujintao.3014@jd.com}

\author[inst1]{Yingzheng Ma}
\ead{mayingzheng.1@jd.com}

\author[inst1]{Jiong Dong}
\ead{dongjiong.1@jd.com}

\author[inst1]{Yongzhi Qi\corref{cor1}}
\ead{qiyongzhi1@jd.com}

\author[inst1]{Jianshen Zhang}
\ead{zhangjianshen@jd.com}

\author[inst1]{Dongyang Geng}
\ead{gengdongyang@jd.com}

\author[inst1]{Anni Zhang}
\ead{zhanganni3@jd.com}

\cortext[cor1]{Corresponding author.}

\affiliation[inst1]{
  organization={Supply Chain Tech Team Y, JD.com},
  op={}
}

\begin{abstract}
Balance-oriented multi-warehouse inventory allocation is a recurring decision problem in large-scale e-commerce supply chains, in which a fixed replenishment quantity is distributed across warehouses to balance post-allocation inventory coverage while accounting for demand forecasts and heterogeneous allocation constraints.
In practice, allocation requirements are often scenario-dependent and expressed in semi-structured or natural-language form rather than as ready-to-solve operations research (OR) formulations. 
We propose an OR-guided Large Language Model (LLM) for Allocation (ORLA) that uses solver feedback to generate, verify, and select OR formulations.
ORLA integrates automatic ``Problem--Model--Code (PMC)'' generation, learning-based formulation selection, and feasibility restoration. We develop three complementary mixed-integer programming formulation families based on deviation minimization, soft band compliance, and knapsack-inspired allocation, together with solver-ready mixed-integer linear programming reformulations, modular constraint extensions, and a penalty-based relaxation mechanism for infeasible cases. The LLM component generates candidate formulations and executable solver code from textual or semi-structured specifications, while the solver provides verification signals for executability, feasibility, and solution quality. To address instance heterogeneity, ORLA estimates the expected quality of candidate formulations, selects promising candidates, and combines their outputs through score-aware aggregation. 
Experimental results on 29 production evaluation batches from JD.com show that the best single OR formulation improves allocation accuracy by 3.4 percentage points over the incumbent approach, while the full ORLA framework achieves a 4.5 percentage-point overall improvement and improves allocation accuracy in 26 of the 29 evaluation batches.
\end{abstract}


\end{frontmatter}

\section{Introduction}

Multi-warehouse inventory allocation is a fundamental decision problem in operations management, supply chain planning, and  operations research (OR). Given a fixed replenishment quantity, current on-hand inventory, forecast demand, and a collection of domain-specific allocation rules, the decision maker must determine how many units to allocate to each warehouse. This problem is closely related to classical topics such as inventory control, replenishment planning, and supply chain operations \citep{Silver1998,Zipkin2000,Axsater2015,HOWARD2011298,SimchiLevi2008}. 

In this paper, a balance-oriented multi-warehouse inventory allocation problem under heterogeneous problem-specific side constraints is studied. The central modeling objective is to allocate a fixed replenishment quantity across warehouses so that post-allocation inventory coverage remains balanced, while hard constraints are respected whenever feasible. We measure inventory coverage through \emph{Target Inventory Days (TID)}, defined as the ratio between available inventory and forecast demand, and use TID-based balance as the organizing principle of the optimization problem.

From a modeling perspective, different allocation scenarios motivate different optimization formulation views.
Some scenarios are naturally expressed by minimizing global deviations from a common coverage target. Others are better captured through target-band compliance with controlled violations. Still others require selective allocation under tightly coupled budget and service constraints. These differences suggest the need for a family of complementary formulations rather than a single universal model.

Another challenge is that heterogeneous allocation requirements are frequently specified in semi-structured or natural-language form rather than as ready-to-solve optimization models. 
This challenge is related to recent work on using 
large language models (LLMs) to translate natural-language descriptions into mathematical formulations and/or solver-compatible code pipelines \citep{Ramamonjison2022,Astorga2024,Jiang2024,Huang2025,AhmadiTeshnizi2024}. 
In contrast to this general line of work, this paper focuses on multi-warehouse inventory allocation problems.

We propose an OR-guided LLM for Allocation (ORLA), a solver-verified methodology that combines exact mixed-integer programming (MIP) modeling, automated formulation generation, feasibility relaxation, and predictive formulation selection. On the OR side, three complementary model families are developed: a \emph{sum-of-deviations} formulation that directly minimizes aggregate imbalance, a \emph{soft band-constrained formulation} that penalizes violations of target coverage intervals, and a \emph{knapsack-inspired formulation} that captures budget-coupled allocation preferences under selective service logic. 
These formulations are further reformulated as mixed-integer linear programs (MILPs) and extended through a modular library of heterogeneous side constraints. To handle the case in which strict requirements render the OR model infeasible, a penalty-based relaxation mechanism is further introduced.

On the artificial intelligence (AI) side, each decision instance is represented as a Problem--Model--Code (PMC) triple, consisting of a textual problem specification, an optimization model, and executable solver code. 
LLMs are attractive in this setting because instruction-tuned and feedback-aligned models have shown strong instruction-following capabilities \citep{Ouyang2022,Bai2022}.
This representation supports an end-to-end pipeline in which an LLM generates candidate optimization formulations and corresponding code from natural-language or semi-structured allocation requirements. 

After that, the generated code is executed to call a MILP solver, and the resulting status and solution are used to verify executability, feasibility, and solution quality. These outputs are further used to detect invalid formulations, construct solver-grounded training signals, and keep the generated decisions auditable.
To improve robustness, supervised fine-tuning (SFT) is adopted, followed by preference optimization based on binary desirability signals \citep{Ouyang2022,Ethayarajh2024}.

To adapt to instance-level heterogeneity, we propose a formulation-selection framework that predicts the expected quality of candidate LLM experts for each instance. Guided by the predicted scores, the framework selects the most promising candidates, and combines their outputs through score-aware aggregation.

Our AI+OR methodology is evaluated on real-world multi-warehouse inventory allocation instances from JD.com. The computational study examines the reliability of solver-oriented code generation, the behavior of feasibility restoration under conflicting constraints, and the generalization of the automated modeling pipeline. 
Furthermore, on real-world allocation scenarios from JD.com, ORLA improves TID-based allocation accuracy by 4.5 percentage points overall, with gains observed in 26 of the 29 evaluation batches.

The remainder of this paper is organized as follows. In Section \ref{sec:related-work}, we review related work on LLM-driven automated OR modeling, post-training methods for LLMs, and MIP modeling. Section \ref{sec:preliminaries} introduces the balance-oriented multi-warehouse inventory allocation problem and provides an overview of the ORLA framework. Section \ref{dec:mip-model} presents the three OR formulation families. Section \ref{sec:pmc-llm} describes the PMC generation pipeline and the formulation-selection mechanism. 
Section \ref{sec:experiments} reports methodological evaluations of code-generation reliability and feasibility relaxation. Section \ref{sec:case-study} presents real-world validation results. Finally, we conclude in Section \ref{sec:conclusion}.

\section{Related Work}
\label{sec:related-work}

\subsection{LLM-Driven OR Modeling}

Recent work has begun to systematically study how LLMs can support OR workflows, especially in translating textual problem descriptions into solver-ready optimization artifacts. The NL4Opt competition formalized this direction by evaluating models on entity extraction and intermediate-representation generation for optimization problems from natural-language descriptions \citep{Ramamonjison2022}. Moving beyond problem understanding, several works have investigated end-to-end autoformulation, in which an LLM generates optimization models, and in some cases executable code \citep{Astorga2024,Jiang2024}. ORLM studies customized training pipelines for automated optimization modeling \citep{Huang2025}, while OptiMUS highlights scalable model generation with solver-in-the-loop verification \citep{AhmadiTeshnizi2024}. From the data-generation perspective, MILP-Evolve demonstrates that LLM-guided evolution can be used to generate diverse MILP classes at scale, thereby alleviating dataset scarcity in optimization-oriented model training \citep{Li2025}.

\subsection{Post-Training Methods for LLMs}
SFT is the standard starting point for instruction-following behavior and structured output generation \citep{Ouyang2022}. Beyond SFT, alignment methods such as reinforcement learning from human feedback (RLHF) \citep{Christiano2017,Ziegler2019,Ouyang2022} and more lightweight offline preference objectives, including direct preference optimization (DPO) \citep{Rafailov2023}, identity preference optimization (IPO) \citep{GheshlaghiAzar2024}, simple preference optimization (SimPO) \citep{Meng2024}, and Kahneman-Tversky optimization (KTO) \citep{Ethayarajh2024}, have been widely used to improve response quality under task-specific desirability criteria.

\subsection{MIP Modeling}
The broad applicability of MIP across production-inventory planning, network design, routing, and facility location has made it one of the dominant modeling paradigms in OR \citep{Eskandarpour2015,Baldacci2012,Melo2009,Kasirzadeh2017}. Three ideas are particularly relevant to our formulation design: Big-$M$ modeling for conditional logic \citep{Williams2013}, knapsack-style modeling for capacity-coupled allocation \citep{Kellerer2004}, and equivalent linearization for handling absolute values within MILP models \citep{Williams2013,BertsimasTsitsiklis1997}.

\section{Problem Setup and Methodology Overview}
\label{sec:preliminaries}

\subsection{Problem Setup}
\label{subsec:problem_setup}

A balance-oriented multi-warehouse inventory allocation problem under heterogeneous side constraints is studied in this paper. Consider a set of warehouses indexed by $k\in\{1,\ldots,n\}$. For a given stock keeping unit (SKU), the input consists of: (i) the on-hand inventory $I_k\ge0$ at warehouse $k$, (ii) the forecasted daily demand $D_k>0$, and (iii) the total replenishment quantity $R\in \mathbb{Z}_{\ge 0}$ to be allocated across warehouses.

The decision is an integer allocation plan $x_k\in \mathbb{Z}_{\ge 0}$ for each warehouse $k$, satisfying the conservation constraint
$\sum_{k=1}^{n} x_k = R$.
The goal is to produce a feasible and high-quality allocation plan that balances inventory coverage across warehouses while satisfying heterogeneous allocation constraints.

\begin{figure}[t]
  \centering
\includegraphics[width=1\linewidth]{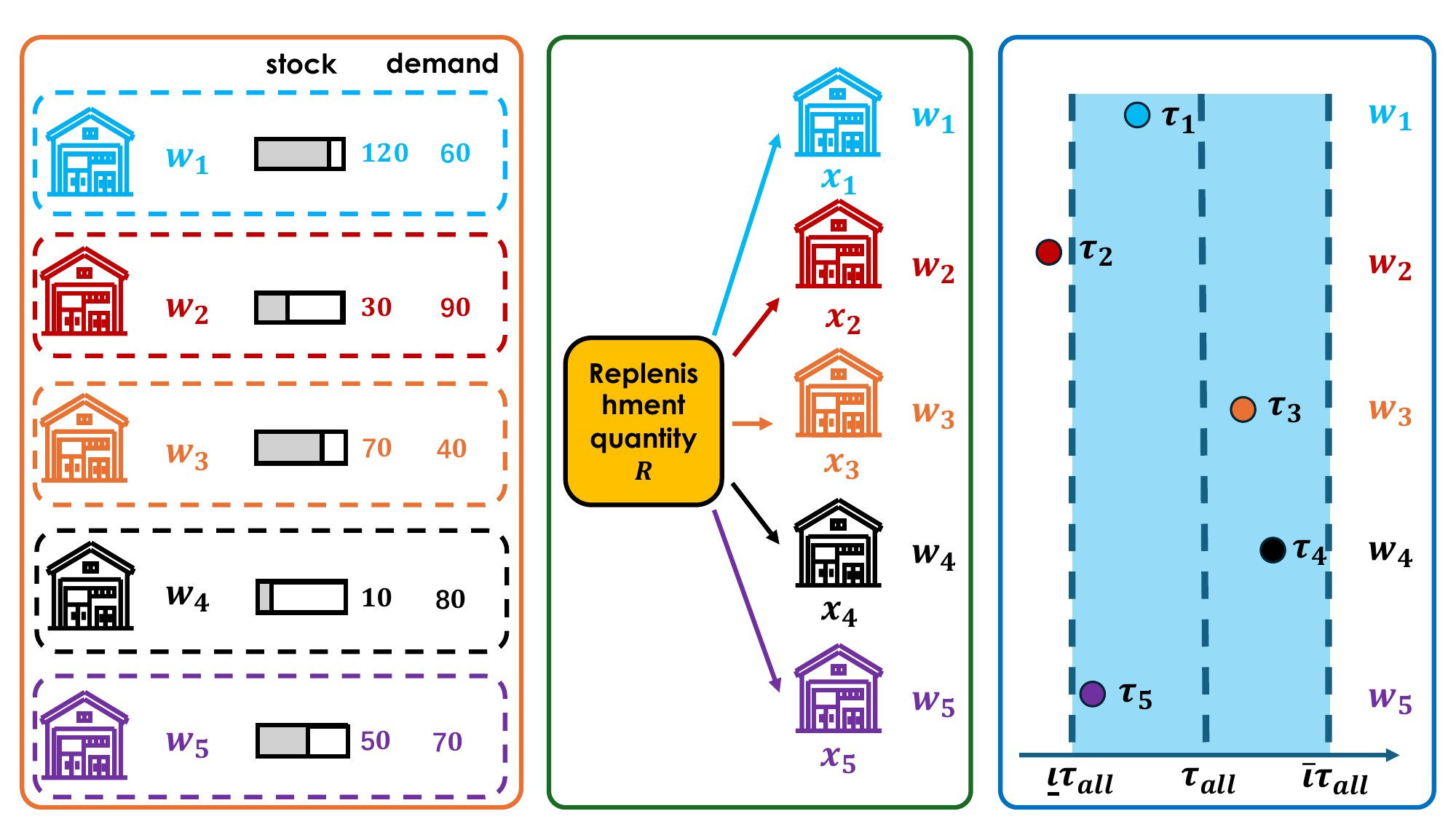}
  \caption{Illustration of TID-based balance in multi-warehouse inventory allocation.}
  \label{fig:DoS}
\end{figure}

For each warehouse $k$, the target inventory days (TID) and target inventory levels (TIL) are defined as
\begin{equation*}
\tau_k := \frac{I_k + x_k}{D_k},~~~
T_k := \tau_{\mathrm{all}}\, D_k,
\end{equation*}
respectively, and the system-level TID is defined as\footnote{Under the replenishment conservation constraint, $\tau_{\rm all}$ and hence $T_{k}$ are constants determined by each instance.}
\begin{equation*}
\tau_{\mathrm{all}} := \frac{\sum_{k=1}^{n} (I_k + x_k)}{\sum_{k=1}^{n} D_k}.
\end{equation*}
To measure the balance quality for each allocation plan, we say that warehouse $k$ is \emph{$(\underline{\ell},\overline{\ell})$-accurately allocated} if its TID lies within a relative band around $\tau_{\mathrm{all}}$, i.e.,
\begin{equation*}
\underline{\ell}\,\tau_{\mathrm{all}} \le \tau_k \le \overline{\ell}\,\tau_{\mathrm{all}},
\end{equation*}
where $0< \underline{\ell}\le 1\le \overline{\ell}$ are given control parameters. Based on these notions, the multi-warehouse allocation accuracy (or balance rate) is defined as
\begin{equation*}
\mathrm{Acc}_{\underline{\ell},\,\overline{\ell}}(x_1,\ldots,x_n)
:=
\frac{
\sum_{k:\ \tau_k\in[\underline{\ell}\,\tau_{\mathrm{all}},\,\overline{\ell}\,\tau_{\mathrm{all}}]}
x_k
}{
\sum_{k=1}^{n}x_k
}.
\end{equation*}
In this paper, the objective is to optimize allocation quality under heterogeneous side constraints, with TID-based balance serving as the central organizing principle of the problem, as shown schematically in Figure \ref{fig:DoS}.

\subsection{Methodology Overview}
\begin{figure}[t]
  \centering
\includegraphics[width=1\linewidth]{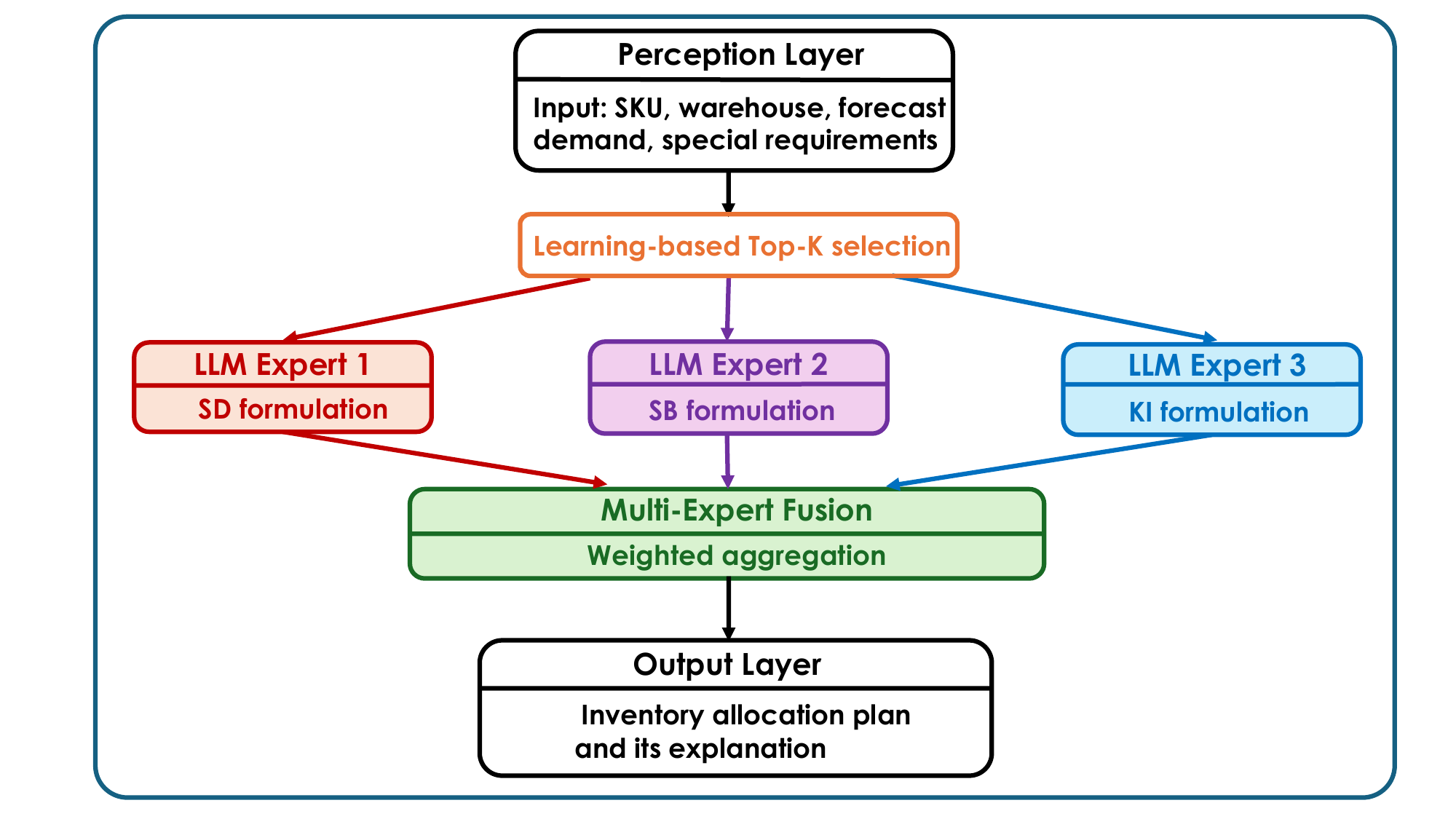}
  \caption{Overview of ORLA framework.}
  \label{fig:system-overview}
\end{figure}

ORLA addresses this problem through a solver-verified methodology that combines exact optimization modeling, automatic formulation generation, learning-based formulation-selection, and feasibility restoration. Figure~\ref{fig:system-overview} summarizes the overall framework.

Accordingly, \emph{three complementary mixed-integer formulation families} are developed in Section~\ref{dec:mip-model}: a sum-of-deviations formulation, a soft band-constrained formulation with controlled slack, and a knapsack-inspired formulation for budget-coupled allocation under selective service logic. Together, these models provide complementary views of the same underlying allocation problem.

To enable automatic generation from textual or semi-structured specifications, we represent each instance as a Problem--Model--Code (PMC) triple and use an LLM to map problem descriptions into candidate \emph{model--code pairs}. The external solver then verifies feasibility and solution quality, serving both as an execution engine and as a source of correctness feedback.

Furthermore, a formulation-selection mechanism is introduced that estimates the quality of candidate LLM experts, selects promising candidates, and combines their outputs.

\section{A Family of OR Formulations}
\label{dec:mip-model}

\begin{figure}[t]
  \centering
\includegraphics[width=1\linewidth]{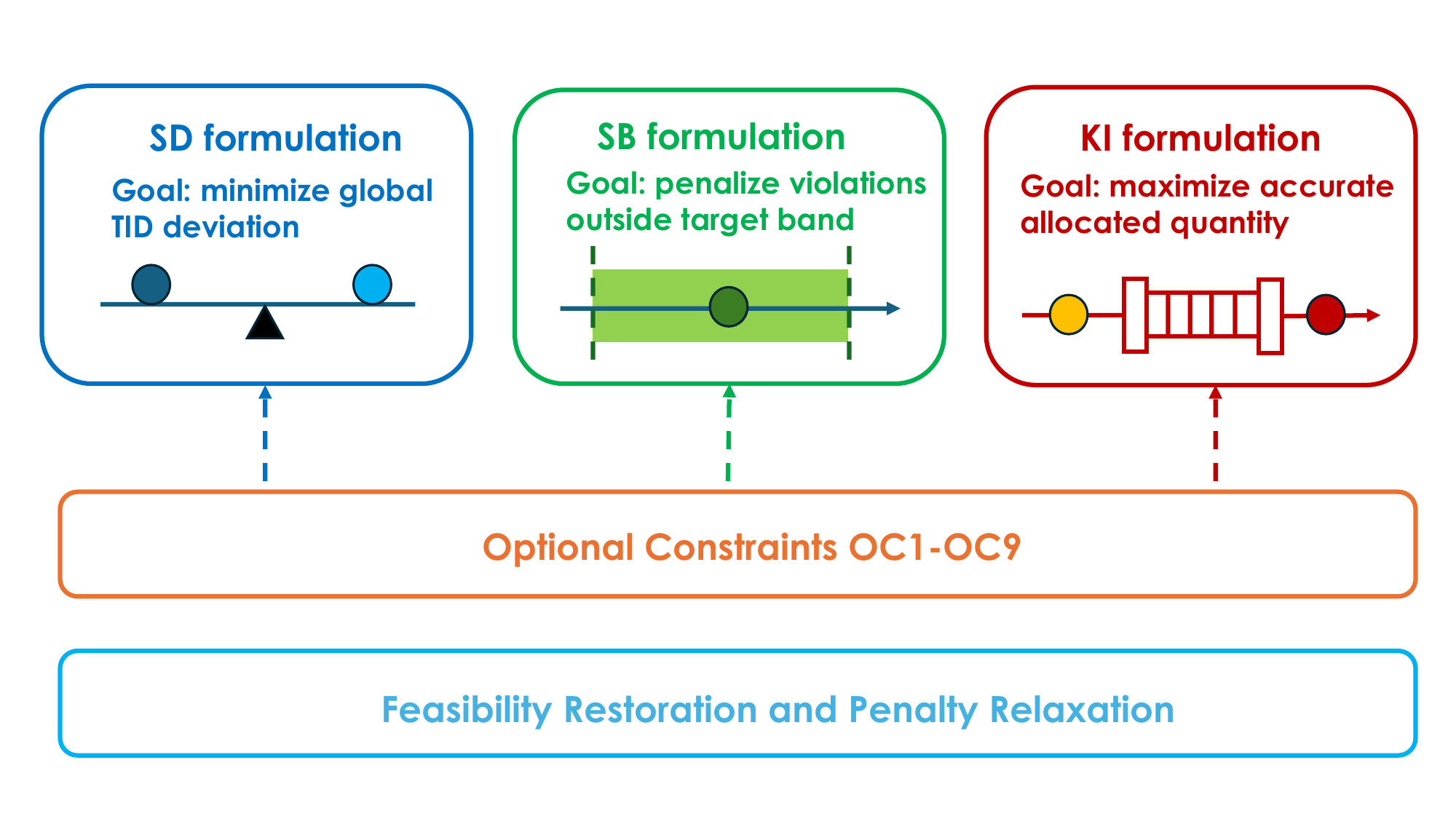}
  \caption{Three complementary formulation families with shared constraints and relaxation modules.}
  \label{fig:mip-families}
\end{figure}

In this section, three complementary mixed-integer formulations are developed, and we further enrich these base models through a modular constraint library and a penalty-based relaxation mechanism for instances in which strict allocation constraints render the original formulation infeasible.
These components form the exact OR backbone of ORLA, as illustrated in Figure \ref{fig:mip-families}.

\subsection{Sum-of-Deviations Formulation}
\label{SDMIP}
\textbf{Basic formulation.}
The \textbf{sum-of-deviations (SD)} MIP problem is formulated as follows.
\begin{equation*}
\begin{aligned}
\min_{\boldsymbol{x}}\quad
& \sum_{k=1}^n
\left|
\frac{I_k+x_k}{D_k} - \tau_{\rm all}
\right| \\
\text{s.t.}\quad
& \sum_{k=1}^n x_k = R, \\
& x_k \in \mathbb{Z}_{\ge 0},\quad k=1, \ldots, n, \\
\end{aligned}
\end{equation*}
where $\boldsymbol{x} = (x_1,
\ldots, x_n)^{\top}$. The objective minimizes the sum of absolute deviations between each warehouse's TID and the global TID benchmark, so that all warehouses stay as close as possible to a common turnover target. 
In addition, the constraints enforce total-quantity conservation and integrality/non-negativity of allocation decisions, ensuring implementable integer shipment plans in practice.

\textbf{Linearized model reformulation.}
The basic formulation above contains absolute-value terms $\lvert\cdot\rvert$ in the objective, which are inconvenient for direct optimization in standard MILP solvers (e.g., SCIP\footnote{https://www.scipopt.org/.}). Therefore, an equivalent linear reformulation is derived as below. 
\begin{equation*}
\begin{aligned}
\min_{{\boldsymbol{x},\,\boldsymbol{\delta}}} \quad
& \sum_{k=1}^n \frac{1}{D_k} \, \delta_k \\
\text{s.t.} \quad
& \delta_k \ge I_k + x_k - T_k, \quad k=1, \ldots, n, \\
& \delta_k \ge -I_k - x_k + T_k, \quad k=1, \ldots, n, \\
& \sum_{k=1}^n x_k = R, \\
& x_k \in \mathbb{Z}_{\ge 0},~~ \delta_k \ge 0, \quad k=1, \ldots, n,
\end{aligned}
\end{equation*}
where $\boldsymbol{x} = (x_1,\ldots, x_n)^\top$, $\boldsymbol{\delta} = (\delta_1,\ldots, \delta_n)^\top$. We refer to this reformulation as the \textbf{linearized sum-of-deviations (LSD)} problem.

\begin{theorem}[\textbf{Equivalence of SD and LSD}]\label{thm:sd-mip}
Assume $D_k>0$ for all $k=1,\ldots,n$. Then SD and its linearized LSD are equivalent.
\end{theorem}

\noindent See \ref{app:proofs-thm1}  for the proof of Theorem~\ref{thm:sd-mip}.
\subsection{Soft Band-Constrained Formulation}
\label{SBMIP}
The core idea of the following \textbf{soft band (SB)} MILP problem is to keep each warehouse's TID within a target interval $[\underline{\ell}\,\tau_{\mathrm{all}},\,\overline{\ell}\,\tau_{\mathrm{all}}]$ while still guaranteeing global replenishment feasibility. 
Instead of enforcing this interval as a hard requirement, nonnegative slack variables are introduced to measure lower-side and upper-side violations, and minimize their total magnitude.

\begin{equation*}
\begin{aligned}
\min_{\boldsymbol{{\rm \xi^+}},\, \boldsymbol{{\rm \xi^-}}, \boldsymbol{x}} \quad 
& \sum_{k=1}^n \xi_k^+ + \sum_{k=1}^n \xi_k^- \\
\text{s.t.} \quad 
& \underline{\ell}\, \tau_{\rm all} - \xi_k^- \le \frac{I_k + x_k}{D_k} \le \overline{\ell}\, \tau_{\rm all} + \xi_k^+, k=1,\ldots,n,\\
& \sum_{k=1}^n x_k = R, \\
& \xi_k^+,\; \xi_k^- \ge 0,\; x_k \in \mathbb{Z}_{\ge 0}, k=1, \ldots, n,
\end{aligned}
\end{equation*}
where $\boldsymbol{\xi^+} = (\xi_1^+, \ldots, \xi_n^+)^\top$, $\boldsymbol{\xi^-} = (\xi_1^-,\ldots,\xi_n^-)^\top$, $\boldsymbol{x} = (x_1, \ldots, x_n)^\top$.

\subsection{Knapsack-Inspired Formulation}
\label{subsec:ki-mip-suite}

\textbf{Basic formulation.}
Inspired by classic knapsack modeling, we incorporate a TID-based balance band through a binary indicator $z_k$ and Big-$M$ constraints, so that the model can explicitly trade off between allocating more units and satisfying balanced allocation requirements. This formulation is referred to as the \textbf{knapsack-inspired (KI)} MIP problem.
\begin{align*}
\max_{\boldsymbol{x},\, \boldsymbol{y}\, ,\boldsymbol{z}} \quad
& \lambda_1 \sum_{k=1}^{n} y_k
- \lambda_2 \sum_{k=1}^{n}
\left| I_k + x_k - T_k \right| \\
\text{s.t.} \quad
& \sum_{k=1}^{n} x_k = R, \\
& I_k + x_k
\ge \underline{\ell} \,T_k - M(1-z_k),
\quad k=1,\ldots,n, \\
& I_k + x_k
\le \overline{\ell} \,T_k + M(1-z_k),
\quad k=1,\ldots,n, \\
& y_k \le x_k,
\quad k=1,\ldots,n, \\
& y_k \le R \, z_k,
\quad k=1,\ldots,n, \\
& y_k \ge x_k - R(1-z_k),
\quad k=1,\ldots,n, \\
& y_k \ge 0,\quad
x_k \in \mathbb{Z}_{\ge 0},\quad
z_k \in \{0,1\},
\quad k=1,\ldots,n .
\end{align*}
where $\boldsymbol{x} = (x_1, \ldots, x_n)^\top$, $\boldsymbol{y} = (y_1, \ldots, y_n)^\top$, $\boldsymbol{z} = (z_1, \ldots, z_n)^\top$.

The second and third constraints activate a TID-based balance requirement when $z_k=1$, forcing warehouse $k$'s post-allocation inventory level to lie within a relative range around the target. When $z_k=0$, the Big-$M$ terms relax this requirement.
Finally, the auxiliary variable $y_k$ is used to reward allocating to warehouses that satisfy the balance band: the lower bound links $y_k$ to $x_k$ when $z_k=1$ and suppresses it when $z_k=0$. A weighted objective is adopted, where the first term maximizes the amount of replenishment assigned to warehouses
whose allocations are counted as band-compliant, whereas the second term
penalizes aggregate deviation from the target inventory levels. 

\textbf{Linearized model reformulation.}
Similarly, nonnegative deviation variables are introduced to obtain a fully solver-ready  \textbf{linearized knapsack-inspired (LKI)} MILP problem as follows. 
\begin{align}
\max_{\boldsymbol{x},\, \boldsymbol{y},\, \boldsymbol{z},\, \boldsymbol{\delta}} \quad
& \lambda_1 \sum_{k=1}^{n} y_k
- \lambda_2 \sum_{k=1}^{n} \delta_k \nonumber\\
\text{s.t.} \quad
& \sum_{k=1}^{n} x_k = R, \tag{C1}\\
& I_k + x_k
\ge \underline{\ell}\, T_k - M(1-z_k),
\quad k=1,\ldots,n, \tag{C2}\\
& I_k + x_k
\le \overline{\ell}\, T_k + M(1-z_k),
\quad k=1,\ldots,n, \tag{C3}\\
& \delta_k
\ge I_k + x_k - T_k,
\quad k=1,\ldots,n, \tag{C4}\\
& \delta_k
\ge T_k - I_k - x_k,
\quad k=1,\ldots,n, \tag{C5}\\
& y_k \le x_k,
\quad k=1,\ldots,n, \nonumber\\
& y_k \le R \, z_k,
\quad k=1,\ldots,n, \nonumber\\
& y_k \ge x_k - R(1-z_k),
\quad k=1,\ldots,n, \nonumber\\
& y_k \ge 0,\quad
\delta_k \ge 0,\quad
x_k \in \mathbb{Z}_{\ge 0},\quad
z_k \in \{0,1\},
\quad k=1,\ldots,n, \nonumber
\end{align}
where $\boldsymbol{\delta} = (\delta_1, \ldots, \delta_n)^\top$.
With similar arguments as in the proof of Theorem \ref{thm:sd-mip}, we know that the above two models are equivalent.
\begin{theorem}[\textbf{Equivalence of KI and linearized LKI}]\label{thm:ki-mip}
Assume $\lambda_1, \lambda_2>0$, and $M\ge \max\limits_{k=1,\ldots,n}\{\underline{\ell}\,T_k - I_k,\, I_k+R-\overline{\ell}\,T_k,\,0\}$. Then KI and the linearized MILP formulation LKI are equivalent.
\end{theorem}
\FloatBarrier

\subsection{Modular Heterogeneous Constraints}

In addition to the balance-oriented requirements captured by the base formulations, real-world decision instances often involve \emph{heterogeneous operational rules}, such as ratio constraints, lower and upper allocation bounds, case-pack restrictions, group-level service requirements, and cardinality controls. Rather than treating these rules as ad hoc modifications for isolated scenarios, we model them as a \emph{modular constraint library} that can be attached systematically to the base formulation families.

The extension patterns of the three formulation families are largely analogous. For clarity, 
we use the KI formulation as a representative example and introduce a collection of optional constraint modules, denoted by \textbf{OC1--OC9} as shown in Table~\ref{tab:ki_mip_modules_old}. 
\begin{table}[t]
\centering
\caption{{Heterogeneous problem-specific side constraints OC1--OC9.}}
\footnotesize
\label{tab:ki_mip_modules_old}
\resizebox{\linewidth}{!}{%
\begin{tabular}{@{}l p{0.31\linewidth} p{0.51\linewidth}@{}}
\toprule
\textbf{Constraint} & \textbf{Meaning} & \textbf{Mathematical form} \\
\midrule
OC1 & Ratio lower-bound & $x_k \ge \alpha_k\,R,\ \forall k\in\mathcal{C}^{\ge}_{\text{ratio}}$ \\
OC2 & Ratio upper-bound & $x_k \le \beta_k \, R, \quad \forall k\in\mathcal{C}^{\le}_{\text{ratio}}$ \\
OC3 & Ratio equality & $x_k = \gamma_k\,R,\ \forall k\in\mathcal{C}^{=}_{\text{ratio}}$ \\
OC4 & Quantity lower-bound & $x_k \ge d_k, \forall k \in \mathcal{C}^{\ge}$ \\
OC5 & Quantity upper-bound & $x_k \le e_k, \forall k \in \mathcal{C}^{\le}$ \\
OC6 & Quantity equality & $x_k = c_k, \forall k \in \mathcal{C}^{=}$ \\
OC7 & Case-pack (size: $p_k\in\mathbb{Z}_{>0}$) & $x_k = p_k\,q_k,q_k \in \mathbb{Z}_{\ge 0}, \quad k=1,\ldots,n$\\
OC8 & Group-level minimum share & $\sum_{l\in \mathcal{S}_k} x_l \ge \eta_k\, R,
\quad k=1,\ldots,K$\\
OC9 & Served warehouse limit & $\sum_{k=1}^{n} s_k \le m, x_k \le Ms_k,  s_k\in\{0,1\}$\\
\bottomrule
\end{tabular}
}
\end{table}
Different combinations of these modules generate a family of extended formulations, see Table~\ref{tab:ki_mip_variants} for details.
\begin{table}[t]
\centering
\small
\caption{Constraint-combination variants built on the base formulation KI. Each variant is obtained by augmenting KI with a subset of additional constraints (OC1--OC9), defined in Table~\ref{tab:ki_mip_modules_old}.}
\label{tab:ki_mip_variants}
\begin{tabular}{@{}l c c c c c c c c c@{}}
\toprule
\textbf{Model} & \textbf{OC1} & \textbf{OC2} & \textbf{OC3} & \textbf{OC4} & \textbf{OC5} & \textbf{OC6} & \textbf{OC7} & \textbf{OC8} & \textbf{OC9} \\
\midrule
KI &  &  &  &  &  &  &  &  &  \\
KIV1 & \cmark &  &  &  &  &  &  &  &  \\
KIV2 &  & \cmark &  &  &  &  &  &  &  \\
KIV3 &  &  & \cmark &  &  &  &  &  &  \\
KIV4 & \cmark & \cmark &  &  &  &  &  &  &  \\
KIV5 & \cmark &  & \cmark &  &  &  &  &  &  \\
KIV6 &  & \cmark & \cmark &  &  &  &  &  &  \\
KIV7 & \cmark & \cmark & \cmark &  &  &  &  &  &  \\
KIV8 &  &  &  & \cmark &  &  &  &  &  \\
KIV9 &  &  &  &  & \cmark &  &  &  &  \\
KIV10 &  &  &  &  &  & \cmark &  &  &  \\
KIV11 &  &  &  & \cmark & \cmark &  &  &  &  \\
KIV12 &  &  &  & \cmark &  & \cmark &  &  &  \\
KIV13 &  &  &  &  & \cmark & \cmark &  &  &  \\
KIV14 &  &  &  & \cmark & \cmark & \cmark &  &  &  \\
KIV15 &  &  &  &  &  &  & \cmark &  &  \\
KIV16 &  &  &  &  &  &  &  & \cmark &  \\
KIV17 &  &  &  &  &  &  &  &  & \cmark \\
\bottomrule
\end{tabular}
\end{table}
As an example, KIV1 is given as follows:
\vspace{-5mm}
\begin{align*}
\max_{\boldsymbol{x},\, \boldsymbol{y}, \, \boldsymbol{z}} \quad
& \lambda_1 \sum_{k=1}^{n} y_k
- \lambda_2 \sum_{k=1}^{n}
\left| I_k + x_k - T_k \right| \\
\text{s.t.} \quad
& \sum_{k=1}^{n} x_k = R, \\
& I_k + x_k
\ge \underline{\ell}\, T_k - M(1-z_k),
\quad k=1,\ldots,n, \\
& I_k + x_k
\le \overline{\ell}\, T_k + M(1-z_k),
\quad k=1,\ldots,n, \\
& y_k \le x_k,
\quad k=1,\ldots,n, \\
& y_k \le R \, z_k,
\quad k=1,\ldots,n, \\
& y_k \ge x_k - R(1-z_k),
\quad k=1,\ldots,n, \\
  & \boldsymbol{x_k \ge \alpha_k\,R,\ \forall k\in\mathcal{C}^{\ge}_{\text{ratio}},}\\
& y_k \ge 0,\quad
x_k \in \mathbb{Z}_{\ge 0},~~z_k \in \{0,1\},~~
 k=1,\ldots,n.
\end{align*}

\subsection{Penalty-Based Relaxation}
\label{sec:feasibility-repair-relaxation-modeling}

A natural consequence of heterogeneous side constraints is that strict formulations may become infeasible. Infeasibility often results from mutually inconsistent lower bounds, equality constraints, or group-share constraints.
To address this issue, we introduce a penalty-based relaxation framework in which selected hard constraints are softened through slack variables and associated violation costs in the objective. This relaxation scheme is formulation-agnostic and can be applied to all model families introduced above. Table~\ref{tab:ki_mip_relaxations}
summarizes the corresponding relaxation patterns\footnote{For maximization OR models,
violation penalties are subtracted from the objective, whereas for minimization
models, the same penalty magnitudes are added to the objective.}.
\begin{table}[H]
\centering
\normalsize
\caption{
Relaxation schemes for constraints, including the corresponding slack-variable modeling and penalty
forms added to the objective.
}
\label{tab:ki_mip_relaxations}
\setlength{\hfuzz}{1pt}

\setlength{\tabcolsep}{5pt}
\renewcommand{\arraystretch}{1.08}
\resizebox{\linewidth}{!}{%
\begin{tabular}{@{} l l c c @{}}
\toprule
    \textbf{Constraints} & \multicolumn{1}{c}{\textbf{Relaxations}} & \textbf{Slack variables} & \textbf{Penalties} \\
\midrule
C2 & $I_k + x_k + u_k^{-} \ge \underline{\ell}\,T_k - M(1-z_k),\ k=1,\ldots,n,\ u_k^-\ge 0$ & $\{u_k^{-}\}_{k=1}^n$ & $-\rho\sum_{k=1}^{n} u_k^{-}$ \\
C3 & $I_k + x_k - u_k^+ \le \overline{\ell}\,T_k + M(1-z_k),\ k=1,\ldots,n,\ u_k^+\ge 0$ & $\{u_k^+\}_{k=1}^n$ & $-\rho\sum_{k=1}^{n} u_k^+$ \\
OC1 & $x_k + u_k^- \ge \alpha_k\,R,\ \forall k\in\mathcal{C}^{\ge}_{\text{ratio}},\ u_k^-\ge 0$ & $\{u_k^-\}_{k\in\mathcal{C}^{\ge}_{\text{ratio}}}$ & $-\rho\sum_{k\in\mathcal{C}^{\ge}_{\text{ratio}}} u_k^-$ \\
OC2 & $x_k - u_k^+ \le \beta_k\,R,\ \forall k\in\mathcal{C}^{\le}_{\text{ratio}},\ u_k^+\ge 0$ & $\{u_k^+\}_{k\in\mathcal{C}^{\le}_{\text{ratio}}}$ & $-\rho\sum_{k\in\mathcal{C}^{\le}_{\text{ratio}}} u_k^+$ \\
OC3 & $\left|x_k-\gamma_k\,R\right|\le u_k,\ \forall k\in\mathcal{C}^{=}_{\text{ratio}},\ u_k\ge 0$ & $\{u_k\}_{k\in\mathcal{C}^{=}_{\text{ratio}}}$ & $-\rho\sum_{k\in\mathcal{C}^{=}_{\text{ratio}}}u_k$ \\
OC4 & $x_k + u_k^- \ge d_k,\ \forall k\in\mathcal{C}^{\ge},\ u_k^-\ge 0$ & $\{u_k^-\}_{k\in\mathcal{C}^{\ge}}$ & $-\rho\sum_{k\in\mathcal{C}^{\ge}} u_k^-$ \\
OC5 & $x_k - u_k^+ \le e_k,\ \forall k\in\mathcal{C}^{\le},\ u_k^+\ge 0$ & $\{u_k^+\}_{k\in\mathcal{C}^{\le}}$ & $-\rho\sum_{k\in\mathcal{C}^{\le}} u_k^+$ \\
OC6 & $\left|x_k-c_k\right|\le u_k,\ \forall k\in\mathcal{C}^{=},\ u_k\ge 0$ & $\{u_k\}_{k\in\mathcal{C}^{=}}$ & $-\rho\sum_{k\in\mathcal{C}^{=}}u_k$ \\
OC7 & $\left|x_k - p_k\,q_k\right| \le u_k$, $q_k\in\mathbb{Z}_{\ge 0}$, $u_k\ge 0$, $k=1,\ldots,n$ & $\{u_k\}_{k=1}^n$ & $-\rho\sum_{k=1}^{n} u_k$ \\
OC8 & 
$\sum_{l\in\mathcal{S}_k}x_l+u_k^-
\ge \eta_k\,R,\ u_k^-\ge 0, k=1,\ldots,K$
& $\{u_k^-\}_{k=1}^K$ & $-\rho \sum_{k=1}^K u_k^-$ \\
OC9 & $\sum_{k=1}^{n} s_k \le m + u^+$, $x_k \le Ms_k$, $s_k\in\{0,1\}$, $u^+\ge 0$ & $u^+$ & $-\rho u^+$ \\
\bottomrule
\end{tabular}%
}
\end{table}

\section{LLM-Driven Formulation Generation and Selection}
\label{sec:pmc-llm}

\subsection{Problem--Model--Code Representation}
\label{subsec:pmc-llm}

Each decision instance is represented as a ``Problem--Model--Code (PMC)'' triple: \emph{a textual problem specification (Problem, P)}, \emph{an exact OR model (Model, M)}, and \emph{executable solver code (Code, C)}. The LLM reads the problem description, generates candidate OR formulations and code, and then relies on an external MILP solver for verification.

\begin{figure}[t]
  \centering
  \includegraphics[width=1\linewidth]{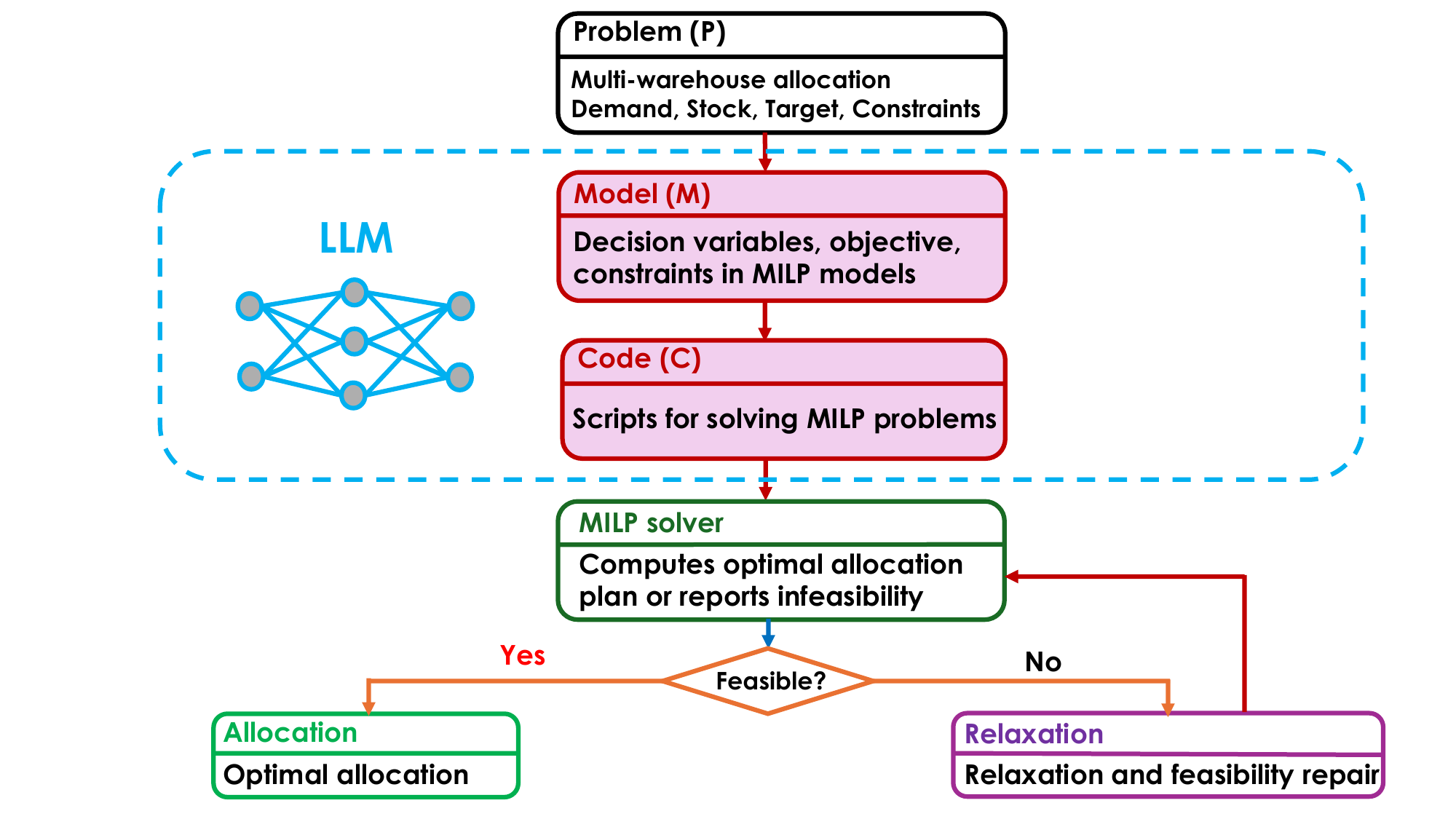}
  \caption{LLM-driven PMC pipeline for multi-warehouse inventory allocation.}
  \label{fig:pmc-example}
\end{figure}

A key design choice is that the solver serves both as an execution engine and as a correctness filter. An \emph{infeasibility-detect--relax--resolve} loop is integrated. The hard-constraint model is solved first. If the solver reports infeasibility, the corresponding relaxation model is activated and re-solved to obtain a feasible allocation plan. This loop is also \emph{central to post-training}: it provides structured supervision for valid formulation/code generation and binary desirability signals for preference optimization. Figure~\ref{fig:pmc-example} illustrates the aforementioned pipeline.

\textbf{Solver-verified supervised fine-tuning.}
Model $\pi_{\rm SFT}$ is obtained through SFT so that it can generate a MIP formulation together with an executable Python script that calls SCIP, from a natural-language multi-warehouse inventory allocation specification. Our SFT dataset contains solver-verifiable PMC triples paired with structured prompts, covering the base formulation and 17 real-world constraint extensions. It also includes positive PMC samples from classical MIP problems such as Max-Cut. Furthermore, we include targeted \emph{negative} instances (e.g., non-linear modeling hallucinations, wrong variable types, missing key constraints) alongside positive instances to improve format robustness and error awareness.

\textbf{Solver-grounded preference optimization.} 
Following the KTO pipeline, a binary-labeled preference dataset from solver verification is constructed. Concretely, we sample a diverse set of multi-warehouse inventory allocation prompts and use the SFT policy $\pi_{\rm SFT}$ to generate a batch of outputs under a high-temperature decoding setting. For each generated candidate, we execute the produced solver script and compare its solution against both (i) the theoretical ground-truth plan derived from the reference MIP and (ii) the required replenishment quantity. 

If the generated code is executable, and its executed allocation matches the ground-truth plan, we label it as a \emph{positive sample}. 
Otherwise, it is labeled as a \emph{negative sample}. Finally, we perform $1\!:\!1$ downsampling between positive and negative samples to obtain a balanced dataset in the standard binary-label format for KTO training.
In practice, it is often difficult to obtain high-quality paired preference samples, while it is relatively easy to accumulate single-label positive or negative samples through continuous annotation and solver-based verification. This data characteristic naturally supports an iterative self-enhancement loop, and is one of the main reasons we adopt KTO as our preference optimization objective.

\subsection{Learning-Based Formulation-Selection}
\label{sec:routing}

A learning-based formulation-selection mechanism over multiple LLM experts is introduced, as shown in Figure~\ref{fig:system-overview}. For each instance, a set of regressors estimates the expected quality of candidate experts, and the most promising Top-$K$ candidates are selected. Their outputs are then combined through score-aware weighting, and feasibility is restored when aggregation introduces mild violations.

\textbf{Predictive estimation.}
For each candidate expert indexed by $i\in\{1,\ldots,E\}$, an
expert-specific LightGBM \citep{lightgbm2017guolin} regressor $\mathrm{LightGBM}_i$ is trained to predict the realized allocation accuracy from instance-level  features such as 
\begin{equation*}
s_{ij} \leftarrow \mathrm{LightGBM}_i(\boldsymbol{\phi}_j), \quad i=1,\ldots,E,
\end{equation*}
where $\boldsymbol{\phi}_j = (\phi_{j_1},\ldots,\phi_{j_l})^\top$ correspond to the features of the inventory allocation instance $j$, and $s_{ij}$ is the predicted score for expert $i$ produced by the trained regressor.

\textbf{Top-$\boldsymbol{K}$ selection and score-weighted combination.}
After generating $E$ scores, we first select the top-$K$ experts. Assume that experts $1,\ldots,K$ are selected. Let $s_{ij}$ denote the predicted score of selected expert $i\in\{1, \ldots, K\}$ for instance $j$. Softmax weights with temperature $\kappa>0$ are computed as below:
\begin{equation*}
\theta_{ij} = \frac{\exp(s_{ij}/\kappa)}{\sum_{\ell=1}^{K}\exp(s_{\ell j}/\kappa)},\qquad i=1,\ldots,K.
\end{equation*}
Let $x_{ij}$ denote the allocation plan produced by selected expert $i$. The fused allocation is computed as
\begin{equation*}
x_j = \sum_{i=1}^{K} \theta_{ij}\,x_{ij}\footnote{
Since the fused solution may be fractional or violate some constraints, we solve a repair MILP problem to restore feasibility while keeping the repaired allocation close to the fused solution.}.
\end{equation*}

\section{Computational Evaluation}
\label{sec:experiments}

This section evaluates the ORLA framework from a methodological perspective. The main questions are: (i) whether post-training improves the reliability of generated code, and (ii) whether the relaxation layer recovers valid solutions under conflicting constraints. The proposed methodology is evaluated on real-world multi-warehouse inventory allocation instances from JD.com supply-chain operations.
For LLM training, we fine-tune the Qwen3-8B\footnote{https://huggingface.co/Qwen/Qwen3-8B.} base model by minimizing the standard token-level cross-entropy loss under teacher forcing, with AdamW \citep{Loshchilov2017} as the optimizer and NEFTune \citep{Jain2023NEFTune} applied during training. The generated MILP problems are solved using SCIP.

\subsection{Data Pools for Post-Training}
\textbf{SFT data pool.} We build PMC-style SFT data for all 18 formulations reported in Table~\ref{tab:ki_mip_variants}, covering the base MIP formulation and its 17 variants. To balance canonical-pattern learning and variant coverage, we set the sampling ratio across positive multi-warehouse inventory allocation samples to $34:1:\cdots:1$. Negative multi-warehouse inventory allocation SFT samples are further constructed by injecting 6 common solver-script error types with approximately uniform frequency across categories. To improve robustness and cross-task generalization, we additionally include PMC samples from 6 classical MIP problems including Max-Cut. In the final mixture, the proportions of multi-warehouse inventory allocation positive samples, multi-warehouse inventory allocation negative samples, and classical MIP samples are 85\%, 10\%, and 5\%, respectively. 

\textbf{KTO data pool.} Starting from $\pi_{\rm SFT}$, we construct binary preference data for offline KTO using solver-verifiable signals. In particular, each candidate is labeled as desirable or undesirable based on executability and solution-quality criteria under the MILP solver, as discussed in Section~\ref{subsec:pmc-llm}.

\subsection{Code Generation Reliability}
\label{subsec:performance-optimization}
In this subsection, two post-training variants are compared, namely SFT-only and SFT+KTO, with a focused evaluation on the code-generation component. Our goal is to quantify how preference optimization affects practical code quality. Specifically, three complementary metrics are reported in Table~\ref{tab:code-quality-sft-kto}: \emph{code syntax accuracy}, \emph{code execution success rate}, and \emph{code execution failure rate}. The evaluation is conducted on a held-out test set containing more than
10,000 real-world multi-warehouse allocation instances.

\begin{table}[H]
  \centering
  \caption{Code quality comparison between SFT and SFT+KTO versions.}
  \label{tab:code-quality-sft-kto}
  \begingroup
  \setlength{\tabcolsep}{10pt}
  \begin{tabular}{lll}
    \toprule
    Evaluation Metric & SFT & SFT+KTO \\
    \midrule
    Code Syntax Accuracy & 99.9\% & \textbf{100\%} \\
    Code Execution Success Rate & 99\% & \textbf{99.9\%} \\
    Code Execution Failure Rate & 1\% & \textbf{0.1\%} \\
    \bottomrule
  \end{tabular}
  \endgroup
\end{table}
On the held-out test set, the SFT+KTO variant attains 100\% code syntax accuracy and reduces the code execution failure rate from 1.0\% to 0.1\% relative to the SFT-only variant.
It is worth noting that these results are obtained under a constrained PMC protocol, where the MILP modeling and solver APIs schema are included during post-training. We observe that preference optimization improves the reliability of the generated solver code on the held-out test set. 
In the following experiments, ``LLM'' refers to the post-trained model obtained via SFT+KTO.

\subsection{Relaxation}
Two multi-warehouse inventory allocation tasks are presented to illustrate how relaxation strategies are applied in practice when hard constraints render the LLM-generated OR model infeasible, as discussed in Section \ref{sec:feasibility-repair-relaxation-modeling}. To fully present the details of model relaxation, we use the linearized model form, i.e., MILP, in this section.

\textbf{Group minimum share task.}
In this task, each predefined warehouse group is required to receive at least a minimum fraction of the total replenishment quantity. Formally, for each group $k$, the model enforces a lower-bound share constraint of the form $\sum_{l\in \mathcal{S}_k} x_l \ge \eta_k\, R,~ k=1,\ldots,K$. Infeasibility can occur when multiple group-level minimum-share constraints are jointly over-restrictive relative to the fixed total budget $R$ and other hard bounds. A typical conflicting case arises when the groups are disjoint and their required minimum shares satisfy $\eta_1+\eta_2>1$.
In our task, we set $R=252$, $(D_1, \ldots, D_8)=(18.1, 9.7, 20.3, 44.7, 17.4, 11.5, 35.2, 8.9)$, $(I_1, \ldots, I_8) = (33, 13, 34, 102, 28, 21, 66, 14)$,  $\eta_1=0.6$, $\eta_2=0.5$, $\underline{\ell}=0.8$, $\overline{\ell}=1.2$, $(\rho_1,\rho_2,\rho_3)=\text{(200, 200, 10000)}$, and $M=1563$. We then relax the corresponding hard constraints and obtain the following relaxation model. 
\begin{tcolorbox}[
enhanced jigsaw,breakable,colback=green!8,colframe=black,boxrule=0.9pt,arc=1mm,left=4pt,right=4pt,top=4pt,bottom=4pt,title={Relaxation model for group minimum share task},colbacktitle=teal!45!black,coltitle=white,fonttitle=\bfseries,title filled]
\footnotesize
\begin{align*}
\max_{\ldots, ~\boldsymbol{u^+, u^-, w^-}}\,
& \lambda_1 \sum_{k=1}^{8} y_k
- \lambda_2\sum_{k=1}^{8} \delta_k
\boldsymbol{- \rho_1 \sum_{k=1}^{8} u_k^+
- \rho_2 \sum_{k=1}^{8} u_k^- - \rho_3 \sum_{k=1}^{2} w_k^-}
\\[-1pt]
\text{s.t.} \quad
& \sum_{k=1}^{8} x_k = R,
\\
& \sum_{l\in \{l_k,l'_k\}}x_l \boldsymbol{+ w_k^-} \ge \eta_k\, R, \quad k = 1,2,\\
& I_k + x_k \boldsymbol{+u_k^-}
\ge \underline{\ell}\,T_k - M(1-z_k),
\quad k = 1,\dots,8,
\\
& I_k + x_k \boldsymbol{- u_k^+}
\le \overline{\ell}\,T_k + M(1-z_k),
\quad k = 1,\dots,8,
\\
& \delta_k \ge I_k + x_k - T_k,
\quad k = 1,\dots,8,
\\
& \delta_k \ge T_k - I_k - x_k,
\quad k = 1,\dots,8,
\\
& y_k \le x_k, y_k \le R z_k,
\quad k = 1,\dots,8,
\\
& y_k \ge x_k - R(1-z_k),
\quad k = 1,\dots,8,
\\
& y_k,\; \delta_k,\; \boldsymbol{u_k^-},\; \boldsymbol{u_k^+ \ge 0}, z_k \in \{0,1\},\; x_k \in \mathbb{Z}_{\ge 0},\quad k = 1,\ldots,8,
\\
& \boldsymbol{w_k^- \ge 0},\, k = 1,2.
\end{align*}
\end{tcolorbox}

\textbf{Allocation quantity constraints task.}
In this task, selected warehouses are required to satisfy explicit quantity constraints, including lower-bound constraints $x_k \ge d_k$ for $k\in\mathcal{C}^{\ge}$ and equality constraints $x_k=c_k$ for $k\in\mathcal{C}^{=}$. Typical conflicting cases include over-constrained equality assignments, incompatible lower bounds across multiple warehouses, or the aggregate required quantity exceeding the available replenishment budget.
The concrete parameter values used in this case are listed as follows: $d_5=45$ for $k\in\mathcal{C}^{\ge}=\{5\}$, $c_6=33$ for $k\in\mathcal{C}^{=}=\{6\}$, $R=72$, $\underline{\ell}=0.8$, $\overline{\ell}=1.2$, and penalty weights $(\rho_1,\rho_2)=(10000, 10000)$. Under the above setting, the original MILP can become infeasible. Then we relax the hard constraints and obtain the following relaxation model.
\begin{tcolorbox}[enhanced,colback=green!8,colframe=black,boxrule=0.9pt,arc=1mm,left=4pt,right=4pt,top=4pt,bottom=4pt,title={Relaxation model for allocation quantity constraints task},colbacktitle=teal!45!black,coltitle=white,fonttitle=\bfseries,title filled]
\footnotesize
\begin{align*}
\max_{\ldots,~\boldsymbol{u^-, v}}\quad
& \lambda_1 \sum_{k=1}^{7} y_k
- \lambda_2\sum_{k=1}^{7} \delta_k
\boldsymbol{- \rho_1 \sum_{k \in \mathcal{C}^{\mathrm{\ge}}} u_k^-
- \rho_2 \sum_{k \in \mathcal{C}^{\mathrm{=}}} v_k}
\\[6pt]
\text{s.t.} \quad
& \sum_{k=1}^{7} x_k = R,
\\[4pt]
& x_k \boldsymbol{+ u_k^-} \ge d_k, \quad k \in \mathcal{C}^{\ge},
\\[4pt]
& x_k - c_k \le \boldsymbol{v_k}, c_k - x_k \le \boldsymbol{v_k}, \quad k \in \mathcal{C}^{=},
\\
& I_k + x_k \ge \underline{\ell}\,T_k - M(1-z_k), \quad k = 1,\dots,7,
\\
& I_k + x_k \le \overline{\ell}\,T_k + M(1-z_k), \quad k = 1,\dots,7,
\\[4pt]
& y_k \le x_k, \quad y_k \le R~ z_k, \quad k = 1,\dots,7,
\\
& y_k \ge x_k - R(1-z_k), \quad k = 1,\dots,7,
\\[4pt]
& \delta_k \ge I_k + x_k - T_k, \quad k = 1,\dots,7,
\\
& \delta_k \ge T_k - I_k - x_k, \quad k = 1,\dots,7,
\\[4pt]
& y_k,\; \delta_k,\; \boldsymbol{u_k^-},\; \boldsymbol{v_k \ge 0}, x_k \in \mathbb{Z}_{\ge 0}, z_k \in \{0,1\}, k=1,\; \ldots,\; 7.
\end{align*}
\end{tcolorbox}

\section{Real-World Validation}
\label{sec:case-study}
We now move from computational evaluation to real-world validation. Section~\ref{subsec:JDaccuracy} reports TID-based allocation accuracy on 29 production evaluation batches from JD.com, comparing the incumbent allocation procedure with the ORLA framework. Section~\ref{sec:JDgeneralization} further demonstrates the generalization capability on task categories not included in post-training, using solver verification of generated code and validity checks of the resulting allocation plans.
\subsection{Production-Batch Allocation Accuracy Evaluation}
\label{subsec:JDaccuracy}

We evaluate the production performance of ORLA on 29 real-world evaluation batches from JD.com. The incumbent allocation procedure is used as the benchmark. We compare it with two OR-driven alternatives. The best-performing single formulation: KI modeling, and the full learning-based predictive formulation-selection method.

\begin{figure}[t]
  \centering
  \includegraphics[width=1\linewidth]{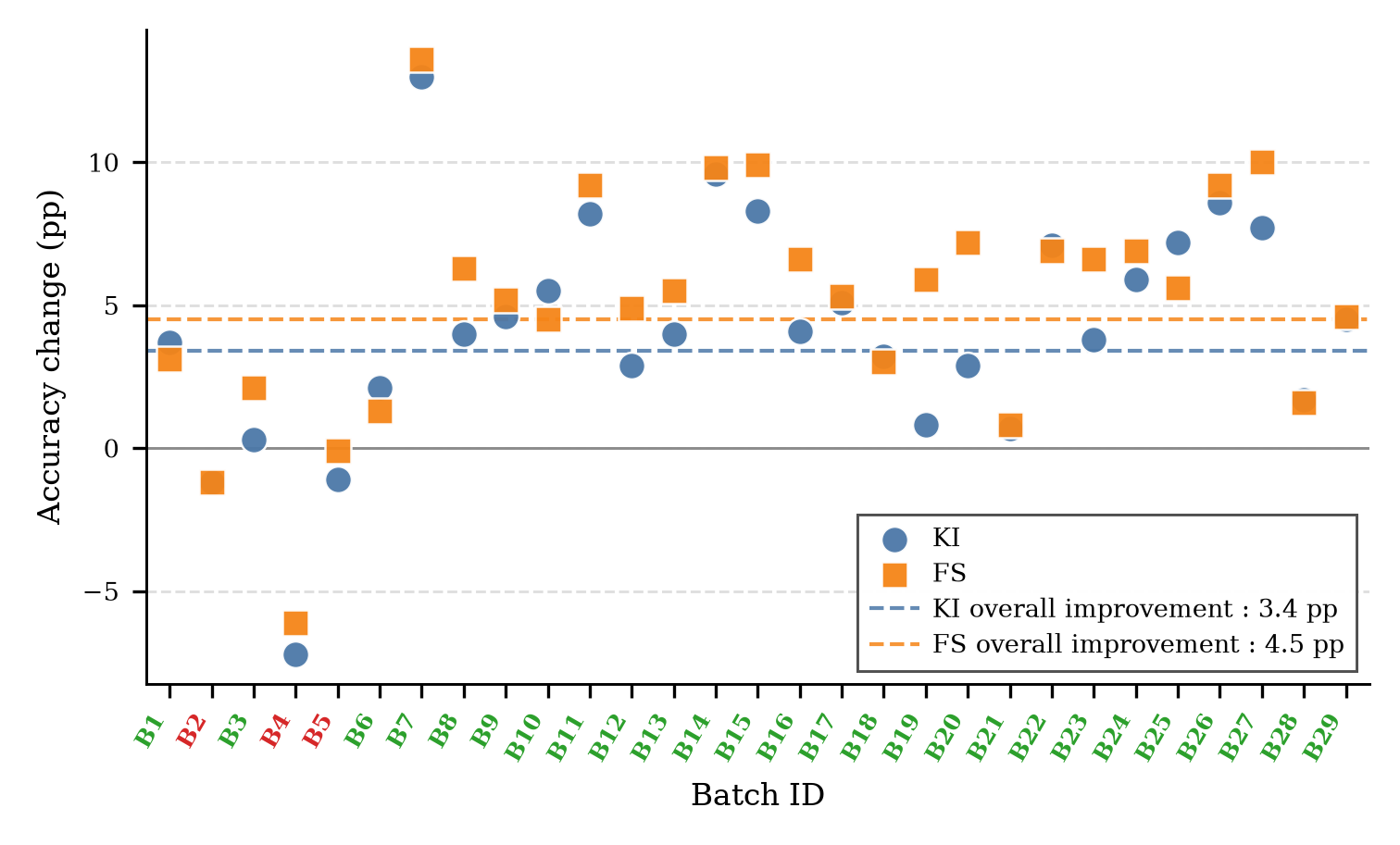}
  \caption{TID-based allocation accuracy changes across 29 evaluation batches. Dashed lines indicate the overall improvements of the KI expert and the formulation-selection (FS) method.}
  \label{fig:accuracy-change-types}
\end{figure}

Figure~\ref{fig:accuracy-change-types} reports the batch-level changes in TID-based allocation accuracy relative to the incumbent. The fixed KI formulation improves 26 of the 29
production batches and achieves an overall improvement of 3.4 percentage points (pp). After applying the learning-based formulation-selection method, the overall improvement increases to 4.5 pp. Thus, both OR-driven approaches improve allocation accuracy in most production batches, and the formulation-selection method provides a
larger aggregate gain than the best fixed formulation.

Detailed batch-level changes are reported in Table~\ref{tab:accuracy-change-26-scenarios}. The formulation-selection method outperforms the fixed KI formulation in 19 of these 26 batches, indicating that its advantage is not driven by a small number of isolated cases but is reflected at the batch level.
A closer inspection shows that the gain from the formulation-selection method is particularly pronounced when the fixed KI formulation provides only moderate improvements. For example, in B19 and B20, the KI formulation improves accuracy by 0.8 and 2.9 pp, respectively, whereas the formulation-selection method increases the corresponding gains to 5.9 and 7.2 pp. This suggests that the formulation-selection method can compensate for cases in which a fixed formulation does not adequately match the structure of a production batch.
The formulation-selection method also preserves and amplifies the high-gain cases of the KI formulation. For all batches in which KI improves accuracy by more than 8 pp, the formulation-selection method achieves larger gains, with additional improvements over KI ranging from 0.2 to 1.6 pp. Moreover, the number of batches with gains
above 5 pp increases from 11 under KI to 17 under the formulation-selection method. These results show that the benefit of the formulation-selection method is not limited to correcting weaker KI cases. It also strengthens performance in batches where the best single formulation is already highly effective.
The downside relative to KI is limited in both frequency and magnitude. KI outperforms the formulation-selection method in 7 of the 26 improved batches, and only two of these cases show a gap of at least one percentage point. Overall, these results show that learning-based formulation-selection increases the upside of the allocation decision while introducing only moderate downside relative to the best fixed formulation.

\begin{table}[h]
  \centering
  \caption{Batch-level allocation accuracy gains of KI and the formulation-selection (FS) method for the 26 production batches with positive gains relative to the incumbent.}
  \label{tab:accuracy-change-26-scenarios}
  \begin{tabular}{lll|lll}
    \toprule
    Batch ID & KI (pp) & FS (pp) & Batch ID & KI (pp) & FS (pp)\\
    \midrule
    B1  & 3.7 & 3.1 & B17 & 5.1 & 5.3 \\
    B3  & 0.3 & 2.1 & B18 & 3.2 & 3.0 \\
    B6  & 2.1 & 1.3 & B19 & 0.8 & 5.9 \\
    B7  & \textbf{13.0} & \textbf{13.6} & B20 & 2.9 & 7.2 \\
    B8  & 4.0 & 6.3 & B21 & 0.7 & 0.8 \\
    B9  & 4.6 & 5.2 & B22 & 7.1 & 6.9 \\
    B10  & 5.5 & 4.5 & B23 & 3.8 & 6.6 \\
    B11  & 8.2 & 9.2 & B24 & 5.9 & 6.9 \\
    B12  & 2.9 & 4.9 & B25 & 7.2 & 5.6 \\
    B13 & 4.0 & 5.5 & B26 & \textbf{8.6} & 9.2 \\
    B14 & \textbf{9.6} & 9.8 & B27 & 7.7 & \textbf{10.0} \\
    B15 & 8.3 & \textbf{9.9} & B28 & 1.7 & 1.6 \\
    B16 & 4.1 & 6.6 & B29 & 4.5 & 4.6 \\
    \bottomrule
  \end{tabular}
\end{table}

\begin{figure}[h]
  \centering
  \includegraphics[width=1\linewidth]{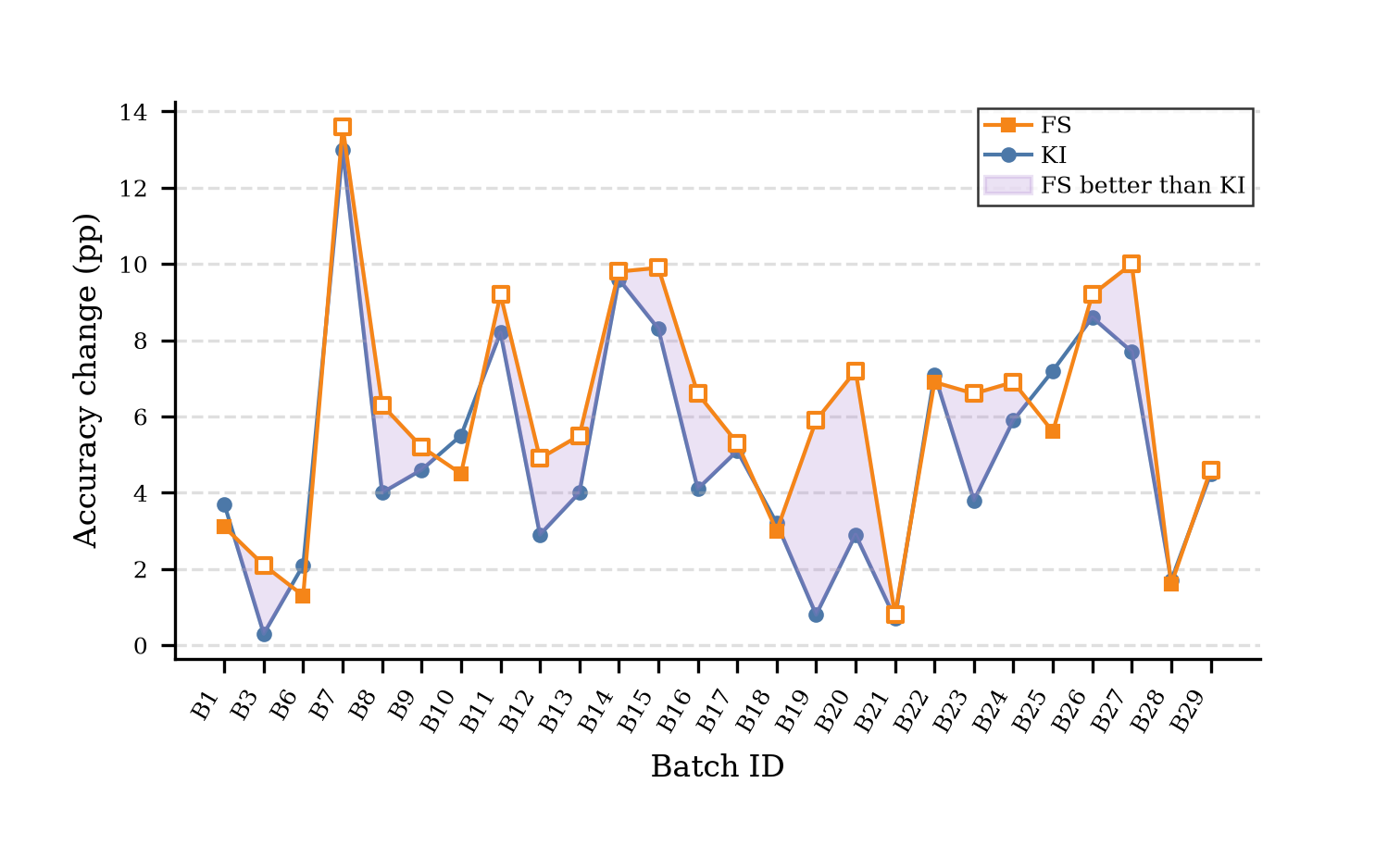}
  \caption{Paired comparison of batch-level allocation accuracy gains between KI and the formulation-selection (FS) method over the 26 improved production batches.
}
  \label{fig:selection-vs-single}
\end{figure}

Figure~\ref{fig:selection-vs-single} presents a paired comparison between the KI formulation and the formulation-selection method over 26 improved batches, focusing on the batch-level response pattern of the formulation-selection method. The two curves exhibit broadly similar movements across production batches, indicating that the formulation-selection method preserves the main allocation behavior induced by the KI formulation. This pattern suggests that the KI formulation captures part of the common structure underlying TID-based accuracy improvements, whereas the formulation-selection method serves as an adaptive refinement mechanism for instance-level heterogeneity. In this sense, the formulation-selection maintains the stable batch-wise behavior of the fixed formulation while selectively modifying allocation decisions when the predicted formulation quality supports such adaptation.

\subsection{Generalization Evaluation}
\label{sec:JDgeneralization}
Beyond the production-batch evaluation, we further evaluate the generalization performance of LLM on unseen task categories. 
For 5 task categories not included in post-training, we run inference with LLM and invoke the solver to verify the generated code, while checking the rationality of the resulting inventory allocation plans.
Table~\ref{tab:generalization} reports results from two dimensions: \emph{code syntax accuracy (CSA)} and \emph{allocation validity (AV)}.
\begin{table}[H]
  \centering
  \caption{Generalization on unseen multi-warehouse inventory allocation tasks}
  \label{tab:generalization}
  \resizebox{\linewidth}{!}{%
  \begin{tabular}{llll}
    \toprule
    Task Name & Task Description & CSA & AV \\
    \midrule
    Equal replenishment & Enforce equal replenishment across warehouses & \cmark & \cmark \\
    Minimum order quantity & Enforce minimum order quantity per warehouse & \cmark & \cmark \\
    Mutually exclusive replenishment & Enforce mutual exclusivity among specified warehouses & \cmark & \cmark \\
    Fixed allocation ratio & Enforce a fixed allocation ratio between two warehouses & \cmark & \xmark \\
    Minimum TID requirement & Enforce minimum TID & \cmark & \xmark \\
    \bottomrule
  \end{tabular}
  }
\end{table}
Our LLM generates syntactically correct and executable solver code for all 5 scenarios. Moreover, for 3 of the 5 scenarios, the generated allocation plans are valid and satisfy all hard constraints. 
The generalization evaluation provides additional evidence beyond the production-batch results. ORLA improves the TID-based allocation accuracy in most real-world evaluation batches, and the LLM component retains partial generalization ability when tested on allocation requirements not included in post-training.

\section{Conclusion}
\label{sec:conclusion}

This paper develops ORLA, a solver-verified and LLM-driven AI+OR methodology for multi-warehouse inventory allocation under heterogeneous domain-specific allocation rules. Three complementary MIP formulations, together with exact MILP reformulations, a modular constraint library, and a formulation-agnostic relaxation scheme, are proposed. Furthermore, we show how solver verification can support structured generation not only at execution time but also during post-training for LLMs through correctness-oriented supervision and preference signals. Empirically,  high reliability in solver-oriented code generation, effective restoration of feasibility under conflicting constraints, and generalization are demonstrated. The real-world case study further shows that the formulation-family approach yields consistent gains over an incumbent allocation process, with additional benefits from instance-dependent formulation selection.

Several directions remain open. First, end-to-end autonomous infeasibility handling can be strengthened
by integrating relaxation and repair strategies directly into post-training, so that the model can proactively diagnose
conflict sources and generate minimal-impact relaxations. Second, the integration of formulation selection logic into model weights
can be explored to enable native formulation-selection capability in LLMs. Third, the framework can be
extended from single-period allocation to multi-period and network-level planning.

\appendix
\section{Proof of Theorem~\ref{thm:sd-mip}}
\label{app:proofs-thm1}

Let
$\Delta_k = I_k+x_k-T_k$  and $\delta_k =\left|\Delta_k \right|, k=1,\ldots,n$
for any feasible allocation $\{x_k\}_{k=1}^n$ of SD problem. 
Then the two inequalities in LSD,
$\delta_k\ge \Delta_k$ and $\delta_k \ge -\Delta_k$
are both satisfied, so $(\{x_k\}_{k=1}^n,\{\delta_k\}_{k=1}^n)$ is feasible for LSD. In this construction, the objective value is exactly the same as that of SD formulation.

Conversely, for any feasible $(\{x_k\}_{k=1}^n,\{\delta_k\}_{k=1}^n)$ of LSD problem, the two inequalities imply $\delta_k\ge \left|\Delta_k\right|$
for each $k$. Since $\delta_k\ge0$, minimizing
$\sum_{k=1}^n \frac{1}{D_k}\delta_k$ forces, at optimum, the smallest feasible value
$\delta_k=\left|\Delta_k\right|$.
Therefore, the two formulations have the same optimal objective value and the same set of optimal allocation solutions.

\section*{Declaration of competing interests}

The authors declare the following financial interests/personal relationships which may be considered as potential competing interests: Yongzhi Qi reports financial support was provided by JD.com Inc. Jintao Xu, Yingzheng Ma, Jiong Dong, Jianshen Zhang, Dongyang Geng, Anni Zhang report financial support was provided by JD.com Inc. The views and opinions expressed in this paper are solely those of the authors and do not necessarily represent those of JD.com Inc.

\bibliographystyle{unsrtnat}
\bibliography{references}

\end{document}